\documentclass[12pt]{article}
\usepackage{amssymb}
\usepackage{amscd}
\usepackage{amsmath}
\usepackage{amsthm}

\let\geq\geqslant
\let\leq\leqslant

\renewcommand{\a}{{\alpha}}
\renewcommand{\b}{{\beta}}
\newcommand{\g}{{\gamma}}

\newcommand{\e}{{\varepsilon}}

\newcommand{\z}{{\zeta}}
\newcommand{\D}{{\Delta}}
\newcommand{\G}{{\Gamma}}

\DeclareMathOperator{\supp}{supp}

\DeclareMathOperator{\Tr}{Tr}

\DeclareMathOperator*{\Res}{Res}

\renewcommand\Im{\hbox{{\rm Im}}\,}
\renewcommand\Re{\hbox{{\rm Re}}\,}
\newcommand{\abs}[1]{\lvert#1\rvert}

\newcommand{\norm}[1]{\lVert#1\rVert}

\newcommand{\R}{{\mathbb R}}
\newcommand{\N}{{\mathbb N}}

\newcommand{\C}{{\mathbb C}}

\numberwithin{equation}{section}

\theoremstyle{plain}
\newtheorem{theorem}{\bf Theorem}[section]
\newtheorem{lemma}[theorem]{\bf Lemma}

\theoremstyle{definition}

\theoremstyle{remark}
\newtheorem*{remark*}{\bf Remark}
\newtheorem{remark}[theorem]{\bf Remark}

\renewcommand{\qed}{\vrule height7pt width5pt depth0pt}

\newcommand{\wt}{\widetilde}

\title{High energy asymptotics and trace formulae for the perturbed harmonic oscillator}

\date{}

\author{Alexander Pushnitski\thanks{Department of Mathematical Sciences,
Loughborough University, Loughborough, LE11 3TU, U.K.
email: a.b.pushnitski@lboro.ac.uk}\mbox{ }
 and Ian Sorrell\thanks{Department of Mathematical Sciences,
Loughborough University, Loughborough, LE11 3TU, U.K.
email: i.sorrell@lboro.ac.uk}
 }

\begin{document}

\maketitle

\begin{abstract}
A one-dimensional quantum harmonic oscillator
perturbed by a smooth compactly supported potential is considered.
For the corresponding eigenvalues $\lambda_n$, a complete
asymptotic expansion for large $n$ is obtained, and the 
coefficients of this expansion are expressed in terms of the 
heat invariants. A sequence of trace formulas is obtained,
expressing regularised sums of integer powers of 
eigenvalues $\lambda_n$ in terms of the heat invariants.
\end{abstract}


\section{Introduction and main results}

\textbf{1. Local heat invariants.}
In order to state our main results, we need to recall the notion of 
\emph{local heat invariants.} Let $v\in C^\infty(\R)$ be any real valued function
such that $v$ and all derivatives of $v$ are uniformly bounded on $\R$.
For the self-adjoint operator  $h=-\frac{d^2}{dx^2}+v$ in  $L^2(\R,dx)$, consider its 
heat kernel $e^{-th}(x,y)$, $t>0$, $x,y\in\R$, i.e. the integral kernel 
of the operator $e^{-th}$.
For any $x\in\R$, one has the asymptotic expansion
\begin{equation}
e^{-th}(x,x)\sim
\frac{1}{\sqrt{4\pi t}}\sum_{j=0}^\infty t^j a_j[v(x)],
\quad t\to+0,
\label{a1}
\end{equation}
where $a_j[v(x)]$ are polynomials in $v$ and derivatives of $v$,
known as the local heat invariants of $h$. 
Explicit formula for $a_j[v(x)]$ is available: 
\begin{equation}
a_j[v(x)]=\sum_{k=0}^{j-1}
\frac{(-1)^j\G(j+\frac{1}{2})}{4^k k! 
(k+j)!(j-k)!\;\G(k+\frac{3}{2})}
(-\tfrac{d^2}{dy^2}+v(y))^{k+j}
(\abs{x-y}^{2k})\mid_{y=x}.
\label{a2}
\end{equation} 
Formula \eqref{a2} was derived in \cite{HitrikPolterovich}
on the basis of the results of \cite{Polterovich1,Polterovich2};
see also references in \cite{Polterovich1} to earlier works on this subject.
From \eqref{a2} or otherwise, one obtains
\begin{align*}
a_0[v(x)]&=1,
\quad
a_1[v(x)]=-v(x),
\quad 
a_2[v(x)]=\frac12 v^2(x)-\frac16 v''(x),
\\
a_3[v(x)]&=-\tfrac16 v^3
+\tfrac16 vv''
+\tfrac1{12}v'^2
-\tfrac1{60}v^{(4)},
\\
a_4[v(x)]&=
\tfrac1{24}v^4
+\tfrac{1}{30}v'v'''
+\tfrac1{60}v v^{(4)}
+\tfrac1{40}(v'')^2
-\tfrac1{840} v^{(6)}
-\tfrac1{12}v''v^2 
-\tfrac1{12}v(v')^2.
\end{align*}

\textbf{2. Perturbed harmonic oscillator.}
Consider the self-adjoint operators
$$
H_0=-\frac{d^2}{dx^2}+x^2
\text{ and }
H=-\frac{d^2}{dx^2}+x^2+q(x)
\text{ in }
L^2(\R,dx),
\text{ where }
q\in C_0^\infty(\R).
$$
These operators can be defined as the closures of the symmetric 
operators, defined on $C_0^\infty(\R)$ by the same differential 
expressions.
Denote by ${\lambda}_n^0=2n-1$, $n=1,2,\dots$ the eigenvalues of $H_0$ 
and by ${\lambda}_1<{\lambda}_2<\cdots$ the eigenvalues of $H$.
The aims of this paper are (i) to describe the asymptotic expansion of 
${\lambda}_n$ as $n\to\infty$, including explicit formulas for the coefficients 
of this expansion in terms of the local heat invariants;
(ii) to derive a series of identities (trace formulas) which relate 
regularized sums of the type $\sum_{n=1}^\infty {\lambda}_n$,  
$\sum_{n=1}^\infty {\lambda}_n^2$, etc.
to some explicit integrals involving heat invariants.

Our results are modelled on the Gel'fand-Levitan-Diki\u\i\mbox{ } trace formulae
for the Sturm-Liouville operator (see \cite{GelfandLevitan,Gelfand,Dikii1,Dikii2}
or \cite{Dikii3}) and in part motivated by the recent advances in
calculation of the heat invariants \cite{Polterovich1,Polterovich2}.

First, as a preliminary result, we establish the asymptotic expansion
\begin{equation}
\Tr(e^{-tH}-e^{-tH_0})
\sim
\frac{1}{\sqrt{4\pi t}}\sum_{j=1}^\infty t^j\int_\R(a_j[x^2+q(x)]-a_j[x^2])dx,
\quad t\to+0,
\label{a4}
\end{equation}
where $a_j$ are the local heat invariants.
In formula \eqref{a4} (as elsewhere in this paper) 
$q\in C_0^\infty(\R)$ and thus all the integrals in the r.h.s. 
converge. 
On the formal level, \eqref{a4} follows by subtracting 
\eqref{a1} with $v(x)=x^2$ from \eqref{a1} with $v(x)=x^2+q(x)$
and integrating over $x$.
A rigourous justification of this formal procedure 
is not difficult and will be given in Section~\ref{sec.c}.

\textbf{3. High energy asymptotics.}
Suppose that $q$ is given. 
Due to the explicit formula \eqref{a2}, we can regard the integrals
appearing in the r.h.s. of \eqref{a4} as known quantities. 
Below we describe the asymptotics of eigenvalues ${\lambda}_n$ in terms of these
integrals. Here is our main result:
\begin{theorem}\label{t1}
(i) One has the asymptotic expansion
\begin{equation}
{\lambda}_n\sim{\lambda}_n^0+\sum_{j=1}^\infty\frac{c_j}{({\lambda}_n^0)^{j/2}},
\quad n\to\infty,
\label{a5}
\end{equation}
with some  coefficients $c_j\in\R$.

(ii) The coefficients $c_j$ in \eqref{a5} can be calculated in the following way.
Consider the formal asymptotic expansion 
\begin{equation}
{\lambda}_n^0\sim{\lambda}_n+\sum_{j=1}^\infty\frac{b_j}{({\lambda}_n)^{j-\frac12}},
\quad n\to\infty,
\label{a6}
\end{equation}
with the coefficients 
\begin{equation}
b_j=(\sqrt{\pi}\G(\tfrac{3}{2}-j))^{-1}\int_\R(a_j[x^2+q(x)]-a_j[x^2])dx.
\label{a7}
\end{equation}
Then inverting the formal asymptotic series \eqref{a6} gives \eqref{a5}.
\end{theorem}

\begin{remark}
1. Theorem~\ref{t1} gives an algorithm of computing the `unknown'
coefficients $c_j$ in the expansion \eqref{a5} in terms of the 
`known' integrals \eqref{a7}. The algorithm is given in the form 
of inverting an asymptotic  series, which is a well defined
algebraic procedure.

In order to compute a coefficient $c_j$, one needs to know 
finitely many coefficients $b_j$.
For example,
$$
c_1=-b_1,\quad
c_2=0,\quad
c_3=-b_2,\quad
c_4=-\frac12 b_1^2,\quad
c_5=-b_3,\quad
c_6=-2b_1 b_2,\quad
c_7=\frac18 b_1^3-b_4.
$$

2. The fact that only half-integer (and not whole integer) 
negative powers of 
${\lambda}_n$ are present in the r.h.s. of \eqref{a6}
is equivalent to a series of identities for the coefficients
$c_j$. For example, the first three identities of this type are
$$
c_2=0, \quad 
c_1^2+2c_4=0, \quad
c_6+c_2^2+2c_1c_3=0.
$$

3. From Theorem~\ref{t1}(ii) we obtain, in particular,
\begin{gather*}
c_1=\frac1\pi \int_\R q(x)dx,\quad c_2=0,\quad
c_3=\frac1\pi\int_\R q(x) x^2dx+\frac{1}{2\pi}\int q^2(x)dx,
\\
c_4=- \frac12 c_1^2, 
\quad
c_5=\frac{1}{16\pi}\int_{-\infty}^\infty (q^3(x)+3q^2(x)x^2+3q(x)x^4+\frac12(q'(x))^2+2q(x))dx.
\end{gather*}
\end{remark}

\textbf{4. Trace formulas.}
As a by-product of our construction, we also obtain  trace formulas 
for the eigenvalues ${\lambda}_n$ and ${\lambda}_n^0$.
This result is a direct analogue of the trace formulas for the Sturm-Liouville
problem due to \cite{GelfandLevitan,Dikii1,Dikii2,Gelfand} and our proof follows 
the reasoning of \cite{Dikii2}.
Let us introduce the Zeta functions
\begin{equation}
Z(s)=\sum_{n=1}^\infty {\lambda}_n^{-s},
\quad 
Z_0(s)=\sum_{n=1}^\infty ({\lambda}^0_n)^{-s},
\quad
\Re s>1.
\label{a8}
\end{equation}
If ${\lambda}_n<0$ for some $n$, then ${\lambda}_n^{-s}$ 
should be understood as 
$\abs{{\lambda}_n}^{-s}e^{-i\pi s}$.
If ${\lambda}_n=0$ for some $n$, then the corresponding term 
in the sum $\sum_{n=1}^\infty {\lambda}_n^{-s}$ is omitted.

Due to the explicit formula ${\lambda}_n^0=2n-1$, we have 
$Z_0(s)=(1-2^{-s})\zeta(s)$, where $\z(s)$ is the 
Riemann Zeta function.
By the properties of $\zeta$, we conclude that 
$Z_0(s)$ has a meromorphic continuation 
into the whole complex plane 
with the only pole at $s=1$, and this pole is simple.
The real zeros of $Z_0$ are at $s=-2n$, $n=0,1,2,\dots$.
\begin{theorem}\label{t2}
The function  $Z(s)$ admits meromorphic 
continuation into the whole complex plane.
Its poles are simple and located at 
 $s=1$ and at $s=-\frac12$, $-\frac32$, $-\frac52$, \dots.
We have the identities:
\begin{equation}
Z(-k)=Z_0(-k), \quad k\in\N.
\label{a9}
\end{equation}
\end{theorem}
As in \cite{Dikii1,Dikii2}, formula 
\eqref{a9} can be combined with the asymptotic 
expansion \eqref{a5} to obtain regularised trace identities
 as follows.
Exponentiating the asymptotics \eqref{a5}, we obtain
for any $\Re s>1$:
\begin{equation}
{\lambda}_n^{-s}\sim
\sum_{j=0}^\infty d_j(s)({\lambda}_n^0)^{-s-(j/2)},
\quad
n\to\infty,
\label{a10}
\end{equation}
where $d_j(s)$ are 
explicit polynomials in $s$ and $c_j$.
 For example,
 \begin{align*}
 d_0(s)&=1,\quad
 d_1(s)=d_2(s)=0,\quad
 d_3(s)=-sc_1,\quad
 d_4(s)=-sc_2,\quad
 d_5(s)=-sc_3,
 \\
 d_6(s)&=-sc_4+\frac{s(s+1)}{2}c_1^2,\quad
d_7(s)=-sc_5+s(s+1)c_1c_2.
\end{align*}
Using this notation, we have for any $k\in\N$:
\begin{equation}
Z(s)=
\sum_{n=0}^\infty 
\{{\lambda}_n^{-s}-\sum_{j=0}^{2k+2} d_j(s)({\lambda}_n^0)^{-s-(j/2)}\}
+
\sum_{j=0}^{2k+2}d_j(s)Z_0(s+(j/2)),
\quad \Re s>1.
\label{a11}
\end{equation}
Now both sides of \eqref{a11} can be meromorphically 
continued into the half-plane $\Re s>-k-\frac12$.
By Theorem~\ref{t2}, the l.h.s. of \eqref{a11}
is analytic at $s=-k$. 
By \eqref{a10}, the same applies to the first term in the 
r.h.s. of \eqref{a11}.
Thus, the second term in the r.h.s. of \eqref{a11}
is also analytic at $s=-k$. 
As $Z_0(s)$ has a pole at $s=1$ (and no other poles), it follows 
that $d_{2k+2}(-k)=0$.
Thus, we obtain 
$$
Z(-k)=
\sum_{n=0}^\infty 
\{{\lambda}_n^{k}-\sum_{j=0}^{2k+1} d_j(-k)({\lambda}_n^0)^{k-(j/2)}\}
+
\sum_{j=0}^{2k+1}d_j(-k)Z_0(-k+(j/2)).
$$
Combined with \eqref{a9}, this yields a series of formulas
\begin{equation}
\sum_{n=0}^\infty 
\{{\lambda}_n^{k}-\sum_{j=0}^{2k+1} d_j(-k)({\lambda}_n^0)^{k-(j/2)}\}
+
\sum_{j=1}^{2k+1}d_j(-k)Z_0(-k+(j/2))=0,
\quad k\in\N.
\label{a12}
\end{equation}
In particular, for $k=1,2,3$ we obtain (taking into account that $Z_0(0)=0$)
\begin{gather}
\sum_{n=0}^\infty ({\lambda}_n-{\lambda}_n^0-\frac{c_1}{\sqrt{{\lambda}_n^0}})+c_1 Z_0(-\tfrac12)
=
0;
\label{a13}
\\
\sum_{n=0}^\infty ({\lambda}_n^2-({\lambda}_n^0)^2
-2c_1\sqrt{{\lambda}_n^0}-\frac{2 c_3}{\sqrt{{\lambda}_n^0}}) 
+2c_1Z_0(-\tfrac12)+2c_3Z_0(\tfrac12)
=0;
\label{a14}
\\
\begin{split}
\sum_{n=0}^\infty 
\left(
{\lambda}_n^3-({\lambda}_n^0)^3-3c_1({\lambda}_n^0)^{3/2}-3c_3({\lambda}_n^0)^{1/2}
-3(c_4+c_1^2)-3c_5 ({\lambda}_n^0)^{-1/2}
\right)& 
\\
+3c_1Z_0(-\tfrac32)+3c_3Z_0(-\tfrac12)
+&3c_5 Z_0(\tfrac12)
=0.
\end{split}
\notag
\end{gather}
Formulas \eqref{a13}, \eqref{a14} (in a slightly different form)
were obtained earlier in \cite{KozlovL}. 
\section{Proof of Theorem~\protect\ref{t1}(ii) }\label{sec.b}
The proof of part (i) of Theorem~\ref{t1} is fairly standard and is based on the 
asymptotic theory of solutions to ODEs and on various explicit 
formulas for parabolic cylinder functions (which give the solutions 
to the ODE corresponding to $q=0$). We give this proof in Sections 5-6.
The proof of part (ii) of Theorem~\ref{t1} is the core of 
our construction and is presented in this 
section. The proof is based on the following 
\begin{lemma}\label{l3}
Let $\lambda_n^0=2 n-1,\ n\in\mathbb{N}$, and let ${\lambda}_1<{\lambda}_2<\cdots$ be a sequence
of real numbers such that $\lambda_n=\lambda_n^0+O(1)$ as $n\to\infty$.
Suppose that an asymptotic expansion
\begin{equation}
\lambda_n^0
\sim
\lambda_n+\sum_{j=1}^\infty 
p_j\lambda_n^{-\alpha_j}
+
\sum_{j=1}^\infty 
q_j\lambda_n^{-j},
\quad n\to\infty,
\label{b1}
\end{equation}
holds true,
where $0\leq\a_1<\a_2<\cdots$ are some non-integer exponents 
and $\{p_j\}\subset\R$, $\{q_j\}\subset\R$.
Then one has the asymptotic expansion
\begin{equation}
\sum_{n=1}^\infty e^{-t\lambda_n}
\sim
\frac1{2t}+
\sum_j \frac{p_j}{2}\Gamma(1-\a_j)t^{\a_j} 
+\frac12\log t\sum_{j=1}^\infty q_j\frac{(-1)^{j}}{(j-1)!}t^{j}
+\sum_{k=1}^\infty r_kt^k
\label{b2}
\end{equation}
as $t\to+0$, with some coefficients $\{r_k\}\subset\R$.
\end{lemma}

\begin{proof}[Proof of Theorem 1.1(ii)]
Given Lemma~\ref{l3} and part (i) of Theorem~\ref{t1},
the proof of Theorem~\ref{t1}(ii) is immediate.
Indeed, inverting the asymptotic expansion \eqref{a5} yields
the expansion of the form
$$
{\lambda}_n^0
\sim 
{\lambda}_n+\sum_{j=1}^\infty b_j {\lambda}_n^{\frac12-j}
+\sum_{j=1}^\infty \wt{b}_j {\lambda}_n^{-j},
\quad 
n\to\infty
$$
with some real coefficients $\{b_j\}$, $\{\wt b_j\}$.
Now using Lemma~\ref{l3} and the explicit formula 
$\sum_{n=1}^\infty e^{-t{\lambda}_n^0}=(2\sinh t)^{-1}$, 
we obtain the asymptotic expansion 
$$
\sum_{n=1}^\infty (e^{-t{\lambda}_n}-e^{-t{\lambda}_n^0})
\sim
\frac1{\sqrt{t}}\sum_{j=1}^\infty\frac{b_j}{2}\G(\tfrac32-j) t^j
+\frac12\log t\sum_{j=1}^\infty \wt{b}_j\frac{(-1)^j}{(j-1)!}t^j
+\sum_{k=1}^\infty \wt r_k t^k
$$
with some real coefficients $\{\wt r_k\}$.
Comparing this to \eqref{a4}, we see that all coefficients 
$\wt{b}_i$ vanish and the coefficients $b_j$ are related to 
the heat invariants by formulas \eqref{a7}.
This completes the proof of Theorem~\ref{t1}(ii).
\end{proof}

In the rest of this section, we prove Lemma~\ref{l3}. 
Broadly speaking, this Lemma can be regarded as a discrete 
analogue of the following version of Watson's Lemma:
\begin{lemma}\label{l4}
Let $\psi:\R\to\R$ be a locally bounded measurable function,
such that $\psi({\lambda})=0$ for all ${\lambda}$ near $-\infty$.
Suppose that $\psi$ has the following asymptotic expansion
\begin{equation}
\psi(\lambda)
=
\sum_j p_j\lambda^{-\alpha_j}+
\sum_{j} q_j\lambda^{-\beta_j}+O(\lambda^{-M})
\quad {\lambda}\to\infty,
\label{b3}
\end{equation}
where $\{\alpha_j\}\subset\R\setminus\N$,
$\{\beta_j\}\subset\N$,
$\{p_j\}\subset\R$, $\{q_j\}\subset\R$
are finite sets and $M>\max(\{\a_j\}\cup\{\b_j\})$, 
$M\in(0,\infty)\setminus\N$.
Then the following asymptotic formula for the Laplace transform
of $\psi$ holds true for $t\to+0$:
\begin{equation}
\int_{-\infty}^{\infty}e^{-t\lambda}\psi({\lambda})d{\lambda}
\sim
\sum_i p_i\Gamma(1-\a_i)t^{\a_i-1}
+(\log t)\sum_{j}
q_j\frac{(-1)^{\b_j}}{(\b_j-1)!}t^{\b_j-1}
+\sum_{0\leq k<M-1}r_k t^{k}+O(t^{M-1})
\label{b4}
\end{equation}
with some coefficients $\{r_k\}$.
\end{lemma}
The proof can be performed, for example, by explicit
computation, checking that each term in the asymptotics  
\eqref{b3} gives the desired contribution to \eqref{b4}.

\begin{proof}[Proof of Lemma~\ref{l3}]

1.
Let 
$$
N({\lambda})=\sharp\{n\mid {\lambda}_n<{\lambda}\},
\quad
N_0({\lambda})=\sharp\{n\mid {\lambda}_n^0<{\lambda}\}.
$$
The main idea of the proof is to approximate $N({\lambda})$
by $N_0(\psi({\lambda}))$, where $\psi$ is a function with the asymptotic
expansion \eqref{b1}.
We construct $\psi$ in terms of its inverse as follows.

The formal inversion of the expansion \eqref{b1}
has the form 
\begin{equation}
{\lambda}_n\sim{\lambda}_n^0+\sum_{j=1}^\infty s_j({\lambda}_n^0)^{-\eta_j},
\quad
n\to\infty,
\label{b5}
\end{equation}
where $0\leq\eta_1<\eta_2<\dots$ and $\{s_j\}\subset\R$.
Fix some sufficiently large $M\in(0,\infty)\setminus\N$; we have 
\begin{equation}
{\lambda}_n={\lambda}_n^0
+
\sum_{\eta_j<M} s_j({\lambda}_n^0)^{-\eta_j}+O(({\lambda}_n^0)^{-M}),
\quad
n\to\infty.
\label{b6}
\end{equation}
Let $\phi\in C^\infty(\R)$ be such that 

(i) $\phi({\lambda})\geq0$ for all ${\lambda}\in\R$ and $\phi({\lambda})=0$ for all ${\lambda}\leq1$;

(ii) $\phi({\lambda})$ is strictly increasing for ${\lambda}>1$;

(iii) $\phi({\lambda})={\lambda}+\sum\limits_{\eta_j<M} s_j {\lambda}^{-\eta_j}$
for all sufficiently large ${\lambda}>0$.

Let $\psi\in C^\infty(0,\infty)$ be such that $\phi(\psi({\lambda}))={\lambda}$
$\forall{\lambda}>0$.
Finally, for ${\lambda}>0$ let us write 
$N_0(\lambda)=\frac12\lambda+{\omega}(\lambda)$, where ${\omega}(\lambda)$ is a 2-periodic function.

With this notation we have:
\begin{multline*}
\sum_{n=1}^\infty e^{-t{\lambda}_n}
=
\int_{-\infty}^\infty e^{-t{\lambda}}dN({\lambda})
=
t\int_{-\infty}^\infty e^{-t{\lambda}}N({\lambda})d{\lambda}
\\
=
t\int_{-\infty}^\infty e^{-t{\lambda}}\left(N({\lambda})-N_0(\psi({\lambda}))\right) d{\lambda}
+
\frac12 t\int_0^\infty e^{-t{\lambda}}\psi({\lambda})d{\lambda}
+
t\int_0^\infty e^{-t{\lambda}}\omega(\psi({\lambda}))d{\lambda}
\\
=: F_1(t)+ F_2(t)+ F_3(t).
\end{multline*}
Below we consider separately the integrals $F_1(t)$, $F_2(t)$ and $F_3(t)$.

2.
Consider $F_2(t)$. By the construction of $\psi$, 
we have the asymptotics
$$
\psi({\lambda})
=
{\lambda}+\sum_{\a_j<M}p_j{\lambda}^{-\a_j}
+
\sum_{j<M}q_j{\lambda}^{-j}+O({\lambda}^{-M})
\quad{\lambda}\to\infty
$$
with the same exponents and coefficients as in \eqref{b1}.
By Lemma~\ref{l4}, we obtain:
\begin{equation}
F_2(t)
\sim
\frac{1}{2t}
+
\sum_{\a_j<M}\frac{p_j}{2}\Gamma(1-\a_j)t^{\a_j}
+
\frac12(\log t)\sum_{j<M}
q_j\frac{(-1)^{j}}{(j-1)!}t^{j}
+
\sum_{0<k<M}r_k t^{k}+O(t^{M}),
\quad t\to+0.
\label{b10}
\end{equation}

3.
Consider $F_1(t)$.
By the construction of $\phi$, we have ${\lambda}_n=\phi({\lambda}_n^0)+O(n^{-M})$,
$n\to\infty$, and so 
$$
F_1(t)
=
\sum_{n=1}^\infty(e^{-t{\lambda}_n}-e^{-t\phi({\lambda}_n^0)})
=
\sum_{n=1}^\infty e^{-t{\lambda}_n}(1-e^{tO(n^{-M})}).
$$

It follows that  $F_1(t)$ has at least $[M]-1$ continuous 
derivatives in $t$ on $[0,\infty)$ and therefore, by the Taylor formula, 
\begin{equation}
F_1(t)
=
\sum_{0\leq k<[M]-1} F_1^{(k)}(0)t^k+o(t^{[M]-1}),
\quad t\to+0.
\label{b8}
\end{equation}

4. 
Let us prove that $F_3$ has continuous derivatives in $t\in[0,\infty)$ of any order,
and so
\begin{equation}
F_3(t)
\sim
\sum_{k=0}^\infty F_3^{(k)}(0) t^k,
\quad t\to+0.
\label{b11}
\end{equation}
Fix $N\in\N$. Integrating by parts $N$ times, we obtain 
\begin{multline*}
F_3(t)
=
t\int_0^\infty e^{-t{\lambda}}\omega(\psi({\lambda}))d{\lambda}
=
t\int_0^\infty e^{-t\phi(\mu)}\phi'(\mu){\omega}(\mu)d\mu
\\
=
-\int_0^\infty (e^{-t\phi(\mu)})'{\omega}(\mu)d\mu
=
(-1)^{N+1}\int_0^\infty (e^{-t\phi(\mu)})^{(N+1)}{\omega}_N(\mu)d\mu,
\end{multline*}
where ${\omega}_N$ is a periodic function.
Using the property (iii) of $\phi$, 
we obtain
$$
(e^{-t\phi(\mu)})^{(N+1)}
=
e^{-t\phi(\mu)}
\{(-t)^{N+1}(\phi'(\mu))^{N+1}
+
\sum_{l=0}^N t^l O({\lambda}^{l-N-2-\eta_1})\}.
$$
It follows that $F_3$ has at least $N-1$ continuous derivatives on $[0,\infty)$.
As $N\in\N$ can be taken arbitrary large, this proves the statement.

5. 
Combining \eqref{b10} -- \eqref{b11},
and using the fact that $M$ can be taken arbitrary large,
we get the desired statement.
\end{proof}

\section{Proof of the asymptotic expansion \protect\eqref{a4}}\label{sec.c}

\textbf{1. Asymptotic expansion \eqref{a1}.}
First let us prove that for any bounded \emph{from below} 
function $v\in C^\infty(\R)$, the asymptotic expansion \eqref{a1},
\eqref{a2} holds true locally uniformly in $x\in\R$.
The expansion \eqref{a1} as such is well known,
but all treatments of this expansion 
in the literature that we are aware of, assume boundedness of $v$,
whereas here we have to deal with potentials of the type $x^2+v(x)$.
Below is a simple argument which shows that the
boundedness from above condition can be lifted. 
Let us fix any $R>0$ and prove that \eqref{a1} holds true uniformly in $x\in[-R,R]$.

Let $\tilde v\in C_0^\infty(\R)$ be such that $\tilde v(x)=v(x)$ for 
all $\abs{x}\leq 4R$, and let $\tilde h=-\frac{d^2}{dx^2}+\tilde v$.
The heat kernel expansion for $C_0^\infty$-potentials is 
certainly well known (see e.g. \cite{HitrikPolterovich}
and references to earlier work therein), and so we have 
$$
e^{-t\tilde h}(x,x)\sim
\frac{1}{\sqrt{4\pi t}}\sum_{j=0}^\infty t^j a_j[v(x)],
\quad t\to+0, \quad \abs{x}\leq R,
$$
uniformly in $x\in[-R,R]$. Thus, it suffices to prove the estimate
\begin{equation}
\sup_{\abs{x}\leq R}\abs{e^{-t\tilde h}(x,x)-e^{-t h}(x,x)}
=
O(e^{-c/t}),
\quad t\to+0,\quad c>0.
\label{c1}
\end{equation}
Let $\chi_R$ be the characteristic function of $(-R,R)$ 
in $\R$, and let $\phi\in C_0^\infty(\R)$ be such that 
$\phi(x)=1$ for $\abs{x}\leq 2R$ and $\phi(x)=0$
for $\abs{x}\geq 3 R$.
Denoting $D=\frac{d}{dx}$, we obtain
\begin{multline*}
\chi_R(e^{-th}-e^{-t\tilde h})\chi_R
=
\chi_R(e^{-th}\phi-\phi e^{-t\tilde h})\chi_R
=
-\int_0^t ds\, \chi_R e^{-(t-s)h}(h\phi-\phi \tilde h)
e^{-s\tilde h}\chi_R
\\
=
\int_0^t ds\, \chi_R e^{-(t-s)h}\phi''
e^{-s\tilde h}\chi_R 
+
\int_0^t ds\, \chi_R e^{-(t-s)h}\phi' D
e^{-s\tilde h}\chi_R 
\\
=
\int_0^t ds\, \chi_R e^{-(t-s)h}\phi''
e^{-s\tilde h}\chi_R 
+
\int_0^t ds\, \chi_R e^{-(t-s)h}\phi' D
e^{s D^2}\chi_R 
\\
-
\int_0^t ds  \int_0^s ds_1\, \chi_R e^{-(t-s)h}\phi' D
e^{(s-s_1) D^2}\tilde v e^{-s_1\tilde  h}\chi_R. 
\end{multline*}
From here, using the explicit formula for the heat kernel
$e^{t\frac{d^2}{dx^2}}$ and the well known estimate
\begin{equation}
\abs{e^{-th}(x,y)}
\leq
\frac{1}{\sqrt{4\pi t}}\exp(-\tfrac{(x-y)^2}{4t}-t\inf_\R v),
\quad 
x,y\in\R,
\label{c2}
\end{equation}
we obtain \eqref{c1}.

\textbf{2. Asymptotic expansion 
\eqref{a4}.} 
Now we are ready to prove the asymptotic expansion 
\eqref{a4}. 
Let $R>0$ be sufficiently large so that $\supp q\subset (-R,R)$.
Let $\chi_{2R}$ be the characteristic function of $(-2R,2R)$ 
and let $\tilde \chi_{2R}=1-\chi_{2R}$.
By the previous step of the proof, it suffices to prove that 
\begin{equation}
\Tr(\tilde \chi_{2R}(e^{-tH}-e^{-tH_0})\tilde \chi_{2R})
=
O(e^{-c/t}),
\quad t\to+0, \quad c>0.
\label{c3}
\end{equation}
By \eqref{c2}, we obtain 
\begin{align*}
\norm{\tilde\chi_{2R} e^{-tH}\chi_R}_{S_2}
&=O(e^{-c/t}),
\quad t\to+0,\quad c>0;
\\
\norm{\chi_R e^{-t H_0}\tilde\chi_{2R}}_{S_2}
&=O(e^{-c/t}),
\quad t\to+0,\quad c>0.
\end{align*}
From these estimates and the formula 
$$
\tilde \chi_{2R}(e^{-tH}-e^{-tH_0})\tilde \chi_{2R}
=
-\int_0^t \tilde \chi_{2R} e^{-(t-s)H}
\chi_R\, q\, \chi_R e^{-s H_0} \tilde \chi_{2R}\, ds
$$
we get the required result \eqref{c3}.

\section{Proof of Theorem~\protect\ref{t2}}

We follow the arguments of \cite{Dikii2}.
First let us assume that ${\lambda}_n\not=0$ for all $n$.
Fix $k\in\N$ and consider formula 
\eqref{a11}.
The second term in the r.h.s. is meromorphic in $\C$
with possible poles at 
$s=1-\frac{j}{2}$, $j=0,1,2,\dots,2k+2$.
The first term in the r.h.s. of \eqref{a11} 
admits analytic continuation into the half-plane
$\Re s>-k-\frac12$.
As $k$ can be taken arbitrary large, it follows 
that $Z$ admits a meromorphic continuation into 
the whole complex plane, all poles of $Z$ are 
simple and located at the points $s=1-\frac{j}{2}$,
$j=0,1,2,\dots$.

Next, from the formula
$$
e^{-t{\lambda}_n}=
\frac1{2\pi i}
\int_{\gamma-i\infty}^{\gamma+i\infty}
(t{\lambda}_n)^{-s} \Gamma(s)ds,
\quad \g>0,\quad (t{\lambda}_n)\in\R\setminus\{0\},
$$
we get 
$$
\sum_{n=1}^\infty (e^{-t{\lambda}_n}-e^{-t{\lambda}_n^0})
=
\frac1{2\pi i}
\int_{\gamma-i\infty}^{\gamma+i\infty}
(Z(s)-Z_0(s))t^{-s}\G(s)ds,
\quad t>0,\quad \g>0.
$$
By a standard argument involving shifting the contour of integration
to the left, the last formula yields the following
asymptotic expansion as $t\to+0$:
\begin{equation}
\sum_{n=1}^\infty (e^{-t{\lambda}_n}-e^{-t{\lambda}_n^0})
\sim
\sum_{j=0}^\infty
\Res_{s=1-(j/2)}\bigl( (Z(s)-Z_0(s))t^{-s}\G(s)\bigr),
\quad t\to+0.
\label{d3}
\end{equation}
First note that $Z(s)-Z_0(s)$ does not have poles at any of the 
points $s=0,-1,-2,\dots$.
Indeed, if $Z(s)-Z_0(s)$ did have a pole at $s=-n$ say, then 
$(Z(s)-Z_0(s))\G(s)$ would have a double pole there and then the 
expansion \eqref{d3} would involve a term $C t^n\log t$.
But by \eqref{a4}, no logarithmic terms actually occur in 
the asymptotic expansion.

Next, by \eqref{a4}, there are no integer 
powers of $t$ in the asymptotic expansion, 
which by the same argument leads to the conclusion that 
$Z(-k)-Z_0(-k)=0$ for all $k=0,1,2,\dots$.

Finally, consider the case when one of the eigenvalues
of $H$ vanishes: ${\lambda}_m=0$.
Then the preceding arguments should be repeated for the 
sequence $\{{\lambda}_n\}$, $n\in\N\setminus\{m\}$.
This leads to the same set of results, apart from the formula
$Z(0)=0$; this should be replaced by $Z(0)=-1$.
\qed

\section{Proof of Theorem~\protect\ref{t1}(i)}

Let us define two solutions $\psi^0_\pm=\psi^0_\pm(x,{\lambda})$
of the equation $-\psi''+x^2\psi={\lambda}\psi$ by 
$$
\psi^0_+(x,{\lambda})=U(-\frac{{\lambda}}2,x\sqrt{2}),
\quad
\psi^0_-(x,{\lambda})=U(-\frac{{\lambda}}2,-x\sqrt{2}),
$$
where $U$ is the parabolic cylinder function (see \cite[\S 19.3]{AS}).
For any $x\in\R$, the solutions $\psi^0_\pm(x,{\lambda})$ are entire functions of ${\lambda}$.
For any ${\lambda}\in\C$, the solutions $\psi_\pm^0(x,{\lambda})$ have the asymptotics 
$$
\psi^0_+(x,{\lambda})=\psi^0_-(-x,{\lambda})=(x\sqrt{2})^{({\lambda}-1)/2}e^{-x^2/2}(1+o(1)),
\quad
x\to+\infty,
$$
and the Wronskian 
$w_0({\lambda})=W(\psi^0_-,\psi^0_+)=(\psi^0_-)'_x\psi^0_+-\psi^0_-(\psi^0_+)'_x$
is given by 
\begin{equation}
w_0({\lambda})
=
\frac{2\sqrt{\pi}}{\G(\frac{1-{\lambda}}{2})}
=
\frac{2}{\sqrt{\pi}}\G(\tfrac{1+{\lambda}}{2})\cos(\tfrac{\pi{\lambda}}{2}).
\label{e1}
\end{equation}
At the eigenvalues ${\lambda}_n^0=2n-1$, the Wronskian $w_0({\lambda})$
vanishes and we have
\begin{equation}
\psi_+^0(x,{\lambda}_n^0)=(-1)^{n+1}\psi^0_-(x,{\lambda}^0_n)=
2^{-(n-1)/2}e^{-x^2/2}H_{n-1}(x),
\label{e1a}
\end{equation}
where $H_n$ is the $n$'th Hermite polynomial.

Next, let $\psi_\pm=\psi_\pm(x,{\lambda})$ be the solutions of the equation
$-\psi''+(x^2+q(x))\psi={\lambda}\psi$, normalised by 
\begin{align*}
\psi_+(x,{\lambda})&=\psi^0_+(x,{\lambda}),\quad x>\sup\supp q,
\\
\psi_-(x,{\lambda})&=\psi^0_-(x,{\lambda}),\quad x<\inf\supp q.
\end{align*}
The eigenvalues ${\lambda}_n$ coincide with the zeros of the Wronskian 
$w({\lambda})=W(\psi_-,\psi_+)$.
In Section~\ref{sec.f} we prove the following Lemma, which describes the 
asymptotics of $w({\lambda})$ as $\Re{\lambda}\to+\infty$. Let $\Omega$ be the half-strip
$$
\Omega=\{{\lambda}\in\C\mid \Re{\lambda}\geq0, \abs{\Im{\lambda}}\leq1\},
$$
for ${\lambda}\in\Omega$ let us denote by $\sqrt{{\lambda}}$ the principal branch of the 
square root, so that $\Re\sqrt{{\lambda}}\geq0$.
\begin{lemma}\label{lma.e.1}
The Wronskian $w({\lambda})$ is analytic in ${\lambda}\in{\Omega}$. The following 
asymptotic expansion holds true:
\begin{equation}
w({\lambda})\sim\frac2{\sqrt{\pi}}\G(\tfrac{1+{\lambda}}{2})
\left(\cos(\tfrac{\pi{\lambda}}{2})\sum_{j=0}^\infty\frac{Q_j}{(\sqrt{{\lambda}})^j}
+
\sin(\tfrac{\pi{\lambda}}{2})\sum_{j=0}^\infty\frac{P_j}{(\sqrt{{\lambda}})^j}\right),
\label{e2}
\end{equation}
as $\abs{{\lambda}}\to\infty$, ${\lambda}\in{\Omega}$.
Here $Q_j,P_j\in\C$ are some coefficients, $Q_0=1$, $P_0=0$.
\end{lemma}
Given Lemma~\ref{lma.e.1}, we can prove Theorem~\ref{t1}(i)
as follows.
Fix any sufficiently small 
 $\e>0$, denote $B_{n,\e}=\{z\mid \abs{z-{\lambda}_n^0}\leq\e\}$,
and let $\G_{n,\e}$ be the contour $\partial B_{n,\e}$ oriented 
anti-clockwise.
By Rouche's Theorem combined with a simple continuity 
argument, we obtain that ${\lambda}_n\in B_{n,\e}$ for all sufficiently 
large $n$. 
Next, the zeros of $w$ in the half-strip ${\Omega}$ coincide with the 
zeros of 
$$
\wt w({\lambda})=\frac{\sqrt{\pi} w({\lambda})}{2\G(\frac{1+{\lambda}}{2})}.
$$
It follows that for all sufficiently large $n$ we have 
\begin{equation}
{\lambda}_n=\frac{1}{2\pi i}\int_{\G_{n,\e}}{\lambda}\frac{\wt w'({\lambda})}{\wt w({\lambda})}d{\lambda}.
\label{e3}
\end{equation}
By analyticity of $w$, the asymptotic expansion \eqref{e2} can be 
differentiated.
Thus, we obtain the following asymptotic 
expansion for ${\lambda}\in\G_{n,\e}$, $n\to\infty$:
\begin{equation}
\frac{\wt w'({\lambda})}{\wt w({\lambda})}
=
-\frac{\pi}{2}\tan\frac{\pi {\lambda}}{2}
+
\frac1{\sqrt{{\lambda}}}g_0({\lambda})
+
\frac1{{\lambda}\sqrt{{\lambda}}}g_1({\lambda})\tan\frac{\pi {\lambda}}{2}
+
\sum_{m=2}^\infty(\tan\frac{\pi{\lambda}}{2})^m(\sqrt{{\lambda}})^{-m+1}g_m({\lambda}),
\label{e4}
\end{equation}
where the functions $g_m({\lambda})$ are analytic in ${\lambda}\in\Omega$ and 
have the asymptotic expansions
\begin{equation}
g_m({\lambda})\sim \sum_{k=0}^\infty c_{mk}(\sqrt{{\lambda}})^{-k},
\quad \abs{{\lambda}}\to\infty, 
\quad {\lambda}\in{\Omega}.
\label{e5}
\end{equation}
Substituting the expansions \eqref{e4} and \eqref{e5} into \eqref{e3}
and computing the integrals of the type 
$\int_{\G_{n,\e}}(\tan\frac{\pi{\lambda}}{2})^m{\lambda}^jd{\lambda}$,
we arrive at the expansion \eqref{a5}.
\qed

\section{Proof of Lemma~\protect\ref{lma.e.1}}\label{sec.f}
Let $x>\sup\supp q$; then 
\begin{equation}
w({\lambda})=(\psi_-(x,{\lambda}))'_x \psi^0_+(x,{\lambda})
-\psi_-(x,{\lambda})(\psi^0_+(x,{\lambda}))'_x.
\label{f1}
\end{equation}
We will use this formula and construct $\psi_-$ in a standard way 
as a solution to the integral equation
\begin{equation}
\psi_-(x,{\lambda})=\psi_-^0(x,{\lambda})+\int_{-\infty}^x G_{\lambda}(x,y)q(y)\psi_-(y,{\lambda})dy,
\label{f2}
\end{equation}
where the integral kernel $G_{\lambda}(x,y)$ is given by 
\begin{equation}
G_{\lambda}(x,y)=-\frac1{w_0({\lambda})}(\psi_+^0(x,{\lambda})\psi_-^0(y,{\lambda})
-\psi_-^0(x,{\lambda})\psi^0_+(y,{\lambda})).
\label{f3}
\end{equation}
The kernel $G_{\lambda}(x,y)$ is an entire function of ${\lambda}$ due to the analyticity 
of $\psi^0_\pm(x,{\lambda})$ and the relation \eqref{e1a}.

Let $R>0$ be sufficiently large so that $\supp q\subset (-R,R)$.
Denote $\Delta=[-2R,2R]$; let $L_{\lambda}:C(\Delta)\to C(\Delta)$ 
be the Volterra type integral operator from \eqref{f2}, 
$$
L_{\lambda}: f(x)\mapsto \int_{-2R}^x G_{\lambda}(x,y)q(y)f(y)dy.
$$
Then the solution of the integral equation \eqref{f2} can be written 
as 
$$
\psi_-=\sum_{n=0}^\infty L_{\lambda}^n\psi^0_-,
$$
and so for the Wronskian \eqref{f1} we have the series representation
$$
w({\lambda})=\sum_{n=0}^\infty W(L_{\lambda}^n\psi^0_-,\psi^0_+)(x),
\quad x\in(R,2R).
$$
\begin{lemma}\label{lma.f.1}
For any $n\in\N$ and any $x\in(R,2R)$, the Wronskian 
$W(L_{\lambda}^n\psi^0_-,\psi^0_+)(x)$ is analytic in ${\lambda}\in{\Omega}$ 
and one has the estimate 
\begin{equation}
\abs{W(L_{\lambda}^n\psi^0_-,\psi^0_+)(x)}
\leq
\frac{C({\lambda})^n}{n!}\abs{\G(\tfrac{1+{\lambda}}{2})},
\quad
C({\lambda})=O(\abs{{\lambda}}^{-1/2}),
\quad\abs{{\lambda}}\to\infty,{\lambda}\in{\Omega}.
\label{f4}
\end{equation}
The asymptotic expansion
\begin{equation}
W(L_{\lambda}^n\psi^0_-,\psi^0_+)(x)
\sim
\G(\tfrac{1+{\lambda}}{2})
\left(\cos(\tfrac{\pi{\lambda}}{2})\sum_{j=n}^\infty\frac{Q_j^{(n)}}{(\sqrt{{\lambda}})^j}
+
\sin(\tfrac{\pi{\lambda}}{2})\sum_{j=n}^\infty\frac{P_j^{(n)}}{(\sqrt{{\lambda}})^j}\right),
\label{f5}
\end{equation}
with some coefficients $Q_j^{(n)}$, $P_j^{(n)}$ holds true as $\abs{{\lambda}}\to\infty$,
${\lambda}\in{\Omega}$.
\end{lemma}
Clearly, Lemma~\ref{lma.e.1} follows from Lemma~\ref{lma.f.1}.

\textit{Proof of Lemma~\ref{lma.f.1}:}
1. It is convenient to introduce two linear combinations $e_+$ and $e_-$ of of the 
solutions $\psi^0_\pm$:
$$
e_+(x,{\lambda})
=
\frac{\sqrt{\pi} 2^{(1-{\lambda})/4}}{\cos(\tfrac{\pi{\lambda}}{2})\G(\tfrac{1+{\lambda}}{4})}
\left(e^{-i\pi({\lambda}+1)/4}\psi^0_+(x,{\lambda})+e^{i\pi({\lambda}+1)/4}\psi^0_-(x,{\lambda})\right),
$$
$e_-(x,{\lambda})=e_+(-x,{\lambda})$. The solutions $e_\pm(x,{\lambda})$ 
are analytic in ${\lambda}\in{\Omega}$
(with removeable singularities at ${\lambda}_n^0$ --- see \eqref{e1a}).
These solutions are chosen so that they satisfy the following 
asymptotic expansions:
\begin{align}
e_\pm(x,{\lambda})&\sim 
e^{\pm i\sqrt{{\lambda}}x}\left(1+\sum_{j=1}^\infty
\frac{R_j^\pm(x)}{(\sqrt{{\lambda}})^j}\right),
\quad{\lambda}\to\infty,\quad {\lambda}\in{\Omega},
\label{f6}
\\
(e_\pm(x,{\lambda}))'_x&\sim 
e^{\pm i\sqrt{{\lambda}}x}\left(\pm i\sqrt{{\lambda}}x+\sum_{j=0}^\infty
\frac{\wt R_j^\pm(x)}{(\sqrt{{\lambda}})^j}\right),
\quad{\lambda}\to\infty,\quad {\lambda}\in{\Omega},
\label{f7}
\end{align}
where $R_j^\pm$, $\wt R_j^\pm$ are polynomials in $x$.
The expansion \eqref{f6} follows directly from 
the formulae 19.9.4, 19.9.5, 19.4.2 of \cite{AS},
and \eqref{f7} is obtained by application of the recurrence
formulas \cite[\S19.6]{AS}.

2.
Let us first prove the bound \eqref{f4}. 
We have 
\begin{equation}
W(L_{\lambda}^n\psi_-^0,\psi_+^0)(x)
=
(L_{\lambda}^n \psi_-^0(x,{\lambda}))'_x
\psi^0_+(x,{\lambda})
-
L_{\lambda}^n\psi^0_-(x,{\lambda})
(\psi^0_+(x,{\lambda}))'_x;
\label{f8}
\end{equation}
let us obtain appropriate bounds for each term 
in the r.h.s. of \eqref{f8}.
Expressing $\psi^0_\pm$ in terms of $e_\pm$, 
\begin{equation}
\psi^0_\pm(x,{\lambda})=\frac{1}{2\sqrt{\pi}i}2^{({\lambda}-1)/4}\G(\tfrac{1+{\lambda}}{4})
\left(e^{i\pi({\lambda}+1)/4}e_\mp(x,{\lambda})-e^{-i\pi({\lambda}+1)/4}e_\pm(x,{\lambda})\right),
\label{f8a}
\end{equation}
and using \eqref{f6}, \eqref{f7}, we obtain
\begin{align}
\norm{\psi_\pm^0(\cdot,{\lambda})}_{C(\D)}&=O(\abs{2^{{\lambda}/4}\G(\tfrac{1+{\lambda}}{4})}),
\quad
\abs{{\lambda}}\to\infty,
\quad {\lambda}\in{\Omega},
\label{f9}
\\
\norm{(\psi_\pm^0(\cdot,{\lambda}))'_x}_{C(\D)}&=
O(\abs{{\lambda}^{1/2}2^{{\lambda}/4}\G(\tfrac{1+{\lambda}}{4})}),
\quad
\abs{{\lambda}}\to\infty,
\quad {\lambda}\in{\Omega}.
\label{f11}
\end{align}
Next, expressing the kernel $G_{\lambda}(x,y)$ in terms of $e_\pm$,
\begin{equation}
G_{\lambda}(x,y)=\frac1{4i}\frac{\G(\tfrac{1+{\lambda}}{4})}{\G(\tfrac{3+{\lambda}}{4})}
(e_+(x,{\lambda})e_-(y,{\lambda})-e_-(x,{\lambda})e_+(y,{\lambda})),
\label{f11a}
\end{equation}
and using the asymptotics \eqref{f6}, we obtain
$$
\sup_{\abs{x}\leq R}\sup_{\abs{y}\leq R}
\abs{G_{\lambda}(x,y)}=O(\abs{{\lambda}}^{-1/2}),
\quad \abs{{\lambda}}\to\infty,
\quad {\lambda}\in{\Omega}.
$$
Using this estimate and the fact that $L_{\lambda}$ is a Volterra type operator,
we obtain
\begin{equation}
\norm{L_{\lambda}^n}_{C(\D)\to C(\D)}
\leq
\frac{C({\lambda})^n}{n!},
\quad 
C({\lambda})=O(\abs{{\lambda}}^{-1/2}),
\quad 
\abs{{\lambda}}\to\infty,
\quad {\lambda}\in{\Omega}.
\label{f12}
\end{equation}
Finally, in order to estimate the term $(L_{\lambda}^n\psi^0_-)'_x$, 
let us introduce the operator $L_{\lambda}':C(\D)\to C(\D)$ by 
$$
L'_{\lambda}: f(x)\mapsto \int_{-R}^x \frac{\partial G_{\lambda}(x,y)}{\partial x}q(y)f(y)dy.
$$
Then $(L_{\lambda}^n\psi^0_-(x,{\lambda}))'_x=L'_{\lambda} L_{\lambda}^{n-1}\psi_-^0(x,{\lambda})$.
Using the asymptotics 
\eqref{f6}, \eqref{f7}, we obtain
\begin{equation}
\norm{L'_{\lambda}}_{C(\D)\to C(\D)}
=O(1),
\quad 
\abs{{\lambda}}\to\infty,
\quad {\lambda}\in{\Omega}.
\label{f13}
\end{equation}
Combining \eqref{f8}, \eqref{f9}--\eqref{f13},
we obtain \eqref{f4}.

3. Let us prove the asymptotic expansion \eqref{f5}.
Using \eqref{f8a}, we obtain
\begin{multline*}
W(L_{\lambda}^n\psi^0_-,\psi^0_+)
=
\frac1{4\pi}2^{({\lambda}-1)/2}\G(\tfrac{1+{\lambda}}{4})^2
\bigl(-W(L_{\lambda}^n e_+,e_+)-W(L_{\lambda}^n e_-,e_-)
\\
+e^{-i\pi(1+{\lambda})/2}W(L_{\lambda}^n e_-,e_+)
+e^{i\pi(1+{\lambda})/2}W(L_{\lambda}^n e_+,e_-)\bigr)
\end{multline*}
Denote
\begin{equation}
g_n^\pm(x,{\lambda})=\frac{L_{\lambda}^n e_\pm(x,{\lambda})}{e_\pm(x,{\lambda})};
\label{f14}
\end{equation}
by \eqref{f6}, the denominator does not vanish 
for all sufficiently large ${\lambda}$.
Using this notation, we obtain
\begin{multline*}
W(L_{\lambda}^n\psi^0_-,\psi^0_+)
=
i\sqrt{\pi}\G(\tfrac{1+{\lambda}}{2})
(e^{i\pi(1+{\lambda})/2}g_n^+-e^{-i\pi(1+{\lambda})/2}g_n^-)
\\
+(g_n^+(x,{\lambda}))'_x O(\abs{\G(\tfrac{1+{\lambda}}{2})})
+(g_n^-(x,{\lambda}))'_x O(\abs{\G(\tfrac{1+{\lambda}}{2})}).
\end{multline*}
It suffices to show that $g_n^\pm$ have the asymptotic expansions
\begin{equation}
g_n^\pm(x,{\lambda})\sim
\sum_{j=n}^\infty \frac{S_j^\pm(x)}{(\sqrt{{\lambda}})^j},
\quad
\abs{{\lambda}}\to\infty,
\quad {\lambda}\in{\Omega}
\label{f15}
\end{equation}
for some coefficients $S_j^\pm\in C(\R)$,
and that for any $x\in(R,2R)$,
\begin{equation}
(g_n^\pm(x,{\lambda}))'_x=O(\abs{{\lambda}}^{-\infty}),
\quad
\abs{{\lambda}}\to\infty,
\quad {\lambda}\in{\Omega}.
\label{f16}
\end{equation}
By the definition of $g_n^\pm$, we have 
\begin{equation}
g_{n+1}^\pm (x,{\lambda})=\int_{-R}^x\frac{G_{\lambda}(x,y)}{e_\pm(x,{\lambda})}
g_n(y,{\lambda})e_\pm(y,{\lambda})q(y)dy.
\label{f17}
\end{equation}
Using this formula, the expression \eqref{f11a} for 
$G_{\lambda}(x,y)$ and the asymptotics \eqref{f6}, the expansion
\eqref{f15} can be easily proven by induction.
The asymptotics \eqref{f16} follows by differentiation of \eqref{f17}.

\end{document}